\documentclass[11pt,reqno]{amsart}

\usepackage{amsgen,amsmath,amstext,amsbsy,amsopn,amsfonts,amssymb, enumerate}

\newcounter{mycounter}  
\newenvironment{noindlist}
 {\begin{list}{\bf Step \arabic{mycounter}.~~}{\usecounter{mycounter} \labelsep=0em \labelwidth=0em \leftmargin=0em \itemindent=0em}}
 {\end{list}}
 
 \newcounter{mycounter2}  
 \newenvironment{noindlistcase}
 {\begin{list}{\bf Case \arabic{mycounter2}.~~}{\usecounter{mycounter2} \labelsep=0em \labelwidth=0em \leftmargin=0em \itemindent=0em}}
 {\end{list}}

\def\char{{\rm char}}
\def\Proof{\noindent{\sl Proof.}\ }
\def\qed{{\hfill $\Box$ \medbreak}}

\newtheorem{defi}{Definition}[section]
\newtheorem{lem}[defi]{Lemma}
\newtheorem{thm}[defi]{Theorem}
\newtheorem{cor}[defi]{Corollary}

\newtheorem{prop}[defi]{Proposition}
\newtheorem{example}[defi]{Example}
\newtheorem*{thmm}{Main Theorem}

\DeclareMathOperator{\la}{\langle}
\DeclareMathOperator{\ra}{\rangle}
\DeclareMathOperator{\su}{\subseteq}
\DeclareMathOperator{\charac}{char}
\DeclareMathOperator{\ad}{ad}

\DeclareMathOperator{\Lie}{{\mathfrak L}}

\DeclareMathOperator{\F}{\mathbb{F}}
\DeclareMathOperator{\bZ}{\mathbb{Z}}

\DeclareMathOperator{\fh}{\mathfrak{H}}
\DeclareMathOperator{\fk}{\mathfrak{k}}
\DeclareMathOperator{\fz}{\mathfrak{Z}}
\DeclareMathOperator{\fg}{\mathfrak{g}}
\DeclareMathOperator{\fN}{\mathfrak{N}}

\DeclareMathOperator{\fm}{\mathfrak{m}}

\begin{document}

\title[Lie solvable enveloping algebras of characteristic two]{Lie solvable enveloping algebras of characteristic two}
\author{\textsc{Salvatore Siciliano}}
\address{Dipartimento di Matematica e Fisica ``Ennio De Giorgi", Universit\`{a} del Salento,
Via Provinciale Lecce--Arnesano, 73100--Lecce, Italy}
\email{salvatore.siciliano@unisalento.it}

\author{\textsc{Hamid Usefi}}
\address{Department of Mathematics and Statistics,
Memorial University of Newfoundland,
St. John's, NL,
Canada, 
A1C 5S7}
\email{usefi@mun.ca}

\begin{abstract} Lie solvable restricted enveloping algebras were characterized by  Riley and  Shalev except when the ground field is of characteristic 2. We resolve the characteristic 2 case here which completes the classification. As an application of our result, we obtain a characterization of ordinary Lie algebras over any field whose enveloping algebra is Lie solvable.

\end{abstract}
\date{\today}

\maketitle{}

\section{Introduction}

Let $A$ be an associative algebra over a field $\F$. Then $A$ can be regarded as a Lie algebra by means of the Lie bracket defined by $[x,y]=xy-yx$, for every $x,y\in A$. The algebra $A$ is said to be Lie solvable if it is solvable as a Lie algebra. 

Lie solvable algebras  have been extensively studied over the years.
There has been a special attention to group algebras. Let $\F G$ be the group algebra of the group $G$ over a field $\F$. Recall that $G$ is said to be $p$-abelian if $p > 0$ and $G^\prime$, the commutator subgroup of $G$, is a finite $p$-group.
Moreover, in the zero characteristic case we say that $G$ is 0-abelian if  it is abelian.
Passi, Passman and Sehgal in \cite{PPS}  proved that a group algebra $\F G$ is Lie solvable if and only if 
either $\char \F \neq 2 $ and $G$ is $p$-abelian or $\char \F=2$ and $G$ has a 2-abelian subgroup of index at most  2. 

Restricted  Lie algebras and $p$-groups enjoy similar properties and so it was of interest to find an analogue of Passi-Passman-Sehgal's result for restricted Lie algebras. Let $L$ be a restricted Lie algebra over a field of positive characteristic $p$ and denote by $u(L)$ the restricted (universal) enveloping algebra of $L$. 
 Riley and  Shalev in early 1990s proved  that  if $p\neq 2$ then $u(L)$ is Lie solvable if and only if $L^\prime$ (the derived subalgebra of $L$) is finite-dimensional and $p$-nilpotent. However, they left out the even characteristic case. The purpose of the present paper is to fill  this gap, thereby completing the classification. Our main result shows that the analogue of group ring case in $p=2$ fails for restricted Lie algebras and indeed,  as we shall see below, the characterizations in $p=2$ case are significantly different.
 
 A polynomial identity (PI) is called non-matrix  if it is not satisfied by the algebra $M_2(\F)$ of 2 by 2 matrices over $\F$.
Note that Lie solvability is a non-matrix PI  provided that $\char \F \neq 2$.
Indeed, if $\char \F=2$  then $M_2(\F)$ is Lie center-by-metabelian.
The non-matrix varieties  of algebras have been extensively studied, see for example \cite{K91, L80, MPR,  R97},
 and enveloping algebras have received special attention in this respect \cite{BRU, BRT, RW99}.
Using the standard PI-theory, like Posner's Theorem, one can deduce that if $R$ is an associative algebra that satisfies a non-matrix PI over a field $\F$ of characteristic $p$   then $[R, R]R$ is nil. If we further assume that $R$ is Lie solvable and $p\neq 2$,  then $[R, R]R$ is nil of bounded index (see \cite{R97}).
Moreover, if we restrict ourselves to $R=u(L)$ then  $R$ satisfies a non-matrix PI if and only if
$[R, R]R$ is nil of bounded index (see \cite{RW99}). However,   if $u(L)$ is Lie solvable and $p=2$ then 
$L^\prime$ may not be even  nil as we shall see below in our main result.

 In order to state the main result,  we recall a few definitions.  
A subset $S$ of $L$ is said to be $p$-nilpotent if there exists $m>0$
 such that $S^{[p]^m}=\{x^{[p]^m}\, \vert \, x\in S\}=0$. We denote by $Z(L)$ the center of $L$. Following \cite{Hoch}, we say that a restricted Lie algebra is strongly abelian if it is abelian and its power mapping is zero.  
In analogy with group rings, we say that  a restricted subalgebra $H$ of $L$ is \emph{p-abelian} if $H^\prime$ is finite-dimensional and $p$-nilpotent. For a subset $X$ of $L$ we denote by $\la X\ra_{\F}$ the vector subspace spanned by $X$. Our main result is the following:

\begin{thmm}
Let  $L$ be a   restricted Lie algebra over a field $\F$ of characteristic 2. Let $\bar{\F}$ be the algebraic closure of $\F$ and set $\Lie=L \otimes_{\F} \bar{\F}$.  Then $u(L)$ is Lie solvable if and only if  
$\Lie$ has  a finite-dimensional 2-nilpotent restricted ideal $I$ such that $\bar \Lie=\Lie/I$ satisfies one of the following conditions:
\begin{enumerate} 
\item[{\normalfont (i)}]  $\bar \Lie$ has an abelian restricted ideal of codimension at most $1$;
\item[{\normalfont (ii)}]  $\bar \Lie$ is nilpotent of class 2 and $\dim \bar \Lie/Z(\bar \Lie)=3$;
\item[{\normalfont (iii)}]  $\bar \Lie=  \la  x_1, x_2, y\ra_{\F} \oplus Z(\bar \Lie)$, where  $[x_1,y]=x_1$,  $[x_2,y]=x_2$, and  $[x_1,x_2]\in Z(\bar \Lie)$;
\item [{\normalfont (iv)}]  $\bar \Lie=  \la  x, y\ra_{\F} \oplus H \oplus Z(\bar \Lie)$, where  $H$ is a strongly abelian finite-dimensional restricted subalgebra of $\bar \Lie$ such that $[x,y]=x$, $[y,h]=h$, and $[x,h]\in Z(\bar \Lie)$  for every $h\in H$;
\item [{\normalfont (v)}] $\bar \Lie= \la  x, y\ra_{\F} \oplus H \oplus Z(\bar \Lie) $, where $H$ is a finite-dimensional abelian subalgebra of $\bar \Lie$ such that $[x,y]=x$, $[y,h]=h$, $[x,h]\in Z(\bar \Lie)$, and $[x,h]^{[2]}=h^{[2]}$, for every $h\in H$.
\end{enumerate}
\end{thmm}

 In Example \ref{examplecod} we show that the extension of the ground field 
 is necessary in the statement of our main theorem. Furthermore, note that the cases (ii)-(v) can occur only when $L^\prime$ is finite-dimensional. In other words, if $u(L)$ is Lie solvable and $L^\prime$ has infinite dimension, then $L$ has a 2-abelian restricted ideal of codimension at most 1.
 
 In the last two decades there has been some interest on the derived length 
of Lie solvable group algebras and enveloping algebras (see \cite{ CSS, RT, Sh1, Sh2, S1,  S2, Sp}), and small characteristics have been considered separately, see for example \cite{KS, R, SS2}. 
It is also worth mentioning that besides the interest on their own, restricted enveloping algebras occur naturally in the study of graded group rings (see e.g. \cite{Q, Sh}). For instance, by using this approach, in \cite{Sh}  Shalev showed that a graded group ring of a finitely generated group ring over a field of characteristic $p>0$ satisfies a polynomial identity if and only if the pro-$p$ completion of $G$ has the structure of a $p$-adic  Lie group. 
 
Finally, let $L$ be a Lie algebra over an arbitrary field $\F$ and let $U(L)$ denote the ordinary universal enveloping algebra of $L$. Necessary and sufficient conditions for $U(L)$ to satisfy a polynomial identity have been found in \cite{B}. Moreover, it is known that if $\F$ has characteristic different from 2, then $U(L)$ is Lie solvable only when $L$ is abelian (see \cite[\S 6, Corollary 6.1]{RS1}). This is no longer true in characteristic 2. As an application of our main theorem,   in the concluding section a description of Lie solvable enveloping algebras in characteristic 2 will be obtained, thereby completing the characterization also in the ordinary case.

\section{Preliminary results and Notation}
An important tool in the proof of our main result is the following theorem, obtained by Passman in \cite{Pa} and   Petrogradsky in \cite{Pe} which characterizes restricted enveloping algebras satisfying a polynomial identity:    
\begin{thm}\label{passman} Let $L$ be a restricted Lie algebra over a field of characteristic $p>0$. Then the restricted enveloping algebra $u(L)$ satisfies a polynomial identity if and only if $L$ possesses a restricted ideal $A$ such that:
\begin{enumerate}
\item [{\normalfont (i)}] $A$ has finite codimension in $L$;
\item [{\normalfont (ii)}] $[A,A]$ is finite dimensional and $p$-nilpotent. 
\end{enumerate}
\end{thm}

The next theorem will also play a crucial role in the sequel. It was proved by  Zalesskii and  Smirnov in \cite{SZ} and, independently, by  Sharma and  Srivastava in \cite{SS}.

\begin{thm}\label{S-Z} Let R be a Lie solvable ring of Lie derived length $t\geq 2$. Then the 
two-sided ideal of R generated by $[[R,R], [R,R]], R]$ is associative
nilpotent of index bounded by a function of $t$.
\end{thm}

Let $S$ be a subset of linear transformations on a finite-dimensional vector space $V$ over a field $\F$. Then $S$ is called \emph{triangularizable} if there exists a chain of $S$-invariant subspaces 
$0=V_0\subseteq V_1\subseteq \ldots V_n=V$
 with $\dim_{\F} V_i=i$ for every $i=0,1,\ldots,n$. In the sequel we will use the following result (see Theorem 1.3.2 of \cite{RR}): 

\begin{thm}\label{triang} Let $L$ be a Lie algebra of linear transformations on a finite-dimensional vector space $V$ over an algebraically closed field. Then $L$ is triangularizable if an only if every element of $L^\prime$ is a nilpotent linear transformation of $V$.   
\end{thm}

Let $L$ be a restricted Lie algebra over a field $\F$ of characteristic $p>0$. For a subset $S$ of $L$ we denote by $\la S\ra_p$
the restricted ideal of $L$ generated by $S$. Moreover, $C_L(S)$ will denote the centralizer of $S$ in $L$. We use the symbols $\zeta_j(L)$ ($j\geq 0$) and $\gamma_i(L)$ ($i\geq 1$), respectively, for the terms of the ascending series and descending series of $L$. An element $x$ of $L$ is said to be \emph{$p$-algebraic} if $\dim_{\F}\la x\ra_p< \infty $; an element which is not $p$-algebraic is called \emph{$p$-transcendental}.   Since we shall deal with the case  $p=2$, our notation adjusts accordingly, that is we shall use the term 2-nilpotent, the symbol $\la S\ra_2$, etc.
Also, longer commutators are  left-normed, that is $[x, y, z]=[[x, y], z]$.

\section{Finite-dimensional derived subalgebra}
In this section we consider the case of restricted Lie algebras having a finite-dimensional derived subalgebra.

\begin{lem}\label{reduced}
Let  $L$ be a finitely generated abelian  restricted Lie algebra over a perfect field $\F$ of  characteristic $p>0$. 
If $L$ is free of nonzero $p$-nilpotent elements then $u(L)$ is a reduced ring.
\end{lem}
\Proof
By \cite[Chapter 4, \S 3, Theorem 3.1]{BMPZ},  $L$ decomposes as 
$$
L=\la x_1\ra_p \oplus \cdots \oplus \la x_h\ra_p \oplus \la y_1\ra_p \oplus \cdots \oplus \la y_k \ra_p, 
$$
where the elements $x_i$ are $p$-transcendental and the elements $y_i$ are $p$-algebraic. Let $H=\bigoplus_{i=1}^h\la x_i \ra_p$ and  $T=\bigoplus_{i=1}^k\la y_i \ra_p$.  Then $u(H)$ is isomorphic to a polynomial $\F$-algebra in $h$ indeterminates. Moreover, as $L$ ha no nonzero $p$-nilpotent elements,  by \cite[Chapter 4, Theorem 4.5.8]{W}  we have that  $T$ is a torus and therefore, by a result due to Hochschild (see \cite{H}),  the algebra $u(T)$ is commutative semisimple. As $u(L)\cong u(H)\otimes_{\F} u(T)$, the claim follows at once. 
\qed

\begin{lem}\label{finite-2}
Let  $L$ be a nilpotent restricted Lie algebra  of class 2 over a perfect field $\F$ of characteristic 2. Suppose that $L^\prime$ is finite-dimensional and not 2-nilpotent. If $u(L)$ is Lie solvable then either $L$ has 
a 2-abelian restricted ideal of codimension 1 or 
$L$ has  a  finite-dimensional 2-nilpotent restricted ideal $I$ such that  $Z(L/I)$ has codimension $3$ in $L/I$.
\end{lem}
\Proof Suppose that  $L$ does not contain any 2-abelian restricted ideal of codimension 1. We proceed by a series of reductive steps:
\begin{noindlist}
\item 
\emph{If $J$ is a finite-dimensional and 2-nilpotent restricted ideal of $L$, then we can replace $L$ with $\bar L=L/J$.} Indeed, the algebra $u(L/J)\cong u(L)/u(L)J$ is Lie solvable,  $\bar L^\prime$ is finite dimensional and not 2-nilpotent, and $\bar L$  does not contain any 2-abelian restricted ideal of codimention 1. Now suppose that we are able to prove that $\bar L$ contains a 2-abelian restricted ideal $\bar I$ such that $Z(\bar{L}/\bar{I})$ has codimension $3$ in $\bar{L}/\bar{I}$. Then we have $\bar I=I/J$ for a suitable restricted ideal of $L$ containing $J$.  Clearly, $I$ is finite-dimensional and 2-nilpotent, and $Z(L/I)$ has codimension $3$ in $L/I$. 
\item 
\emph{ We can assume that $ \langle L^\prime\rangle_2$ is free of nonzero 2-nilpotent elements.}
Let $V$ be the subspace spanned by all $z\in L^\prime$ such that $z$ is $2$-nilpotent. Note that $\langle V \rangle_2$ is a central restricted ideal of $L$. Since $\la V\ra_2$ is finite-dimensional and 2-nilpotent, by the previous step we can replace $L$ with $L/\la V\ra_2$. 

Now, by Theorem \ref{passman} there exists a 2-abelian restricted ideal  of $L$ of finite codimention. Let $A$ be such an ideal of minimal codimension.  Then $\dim L/A\geq 2$. By the previous step we can replace $L$ with $L/\la A^\prime\ra_2$ and thereby assume that $A$ is abelian. Moreover, we have $L^\prime\su Z(L)\su A$. Let $z\in L$. Note that if $[z, A]$ is 2-nilpotent then by the minimality of $\dim L/A$ we must have $z\in A$. We fix  $z\in L\setminus A$ and $y\in A$ such that $[z,y]$ is not 2-nilpotent.

\item\label{[N,L]} \emph{We can assume that   $\dim A/Z(L)= 1$.} Let $N$ be the subspace spanned by all
$a\in A$ such that $[z,a]$ is 2-nilpotent. Note that $\dim A/N$ is finite as $L^\prime$ is finite-dimensional, and $N$ is in fact a restricted ideal of $L$. 
Let $x\in L$ and $b\in N$ and consider 
\begin{align*}
u=[[zby, z], [x, xb], y]=[zb[z, y], x[x,b], y].
\end{align*}
By Theorem \ref{S-Z} the element $u$ is nilpotent and one has $u=[x,b][z,y][[z,x]b+z[b,x], y]=[x,b]^2[z,y]^2$.
Since $[z,y]$ is not nilpotent it follows that $[x, b]$ must be nilpotent. As $\la L^\prime \ra_2$ has no nonzero 2-nilpotent elements, this forces $[x,b]=0$ and so $N$ is central in $L$.  Thus, it will be enough to prove that $\dim A/N=1$. Suppose otherwise and let  $a$ be an arbitrary element in $A$ such that $y$ and $a$ are linearly independent modulo $N$.
 Let $x\in L$ such that $x$ and $z$ are linearly independent modulo $A$.
 Consider the element
\begin{align*}
w=[[ az, y], [xya, x], z]&=[z, y][ [a, x]([x,y][a, z]+[y,z][a,x])].
\end{align*}
Then  $w$ is nilpotent by Theorem \ref{S-Z}, and Lemma \ref{reduced} applied to $\la L^\prime\ra_2$ implies that $w=0$. 
Therefore one has   
$$
([z, y][ a, x])^{2}=[z, y][ a, x][x, y][ a, z].
$$
By symmetry and switching $y$ and $a$, we deduce that 
$$
([z, a][ x, y])^{2}=[z, y][ a, x][x, y][ a, z].
$$
 Thus $([z, y][ a, x])^{2}=([z, a][ x, y])^{2}$ and applying  Lemma \ref{reduced} again implies that  $[z, y][ a, x]=[z, a][ x, y]$. 
 Since $y$ and $a$ are linearly independent modulo $N$, it follows that 
$[z, y]$ and $[z, a]$ are linearly independent. Hence, by the PBW Theorem for restricted Lie algebras (see e.g. \cite[Chapter 2, \S 5, Theorem 5.1]{SF}) there exists $\alpha\in \F$ such that $[x, y]=\alpha[z, y]$ and
$[x, a]=\alpha[z, a]$. Put $x_1=x+\alpha z$. Then one has $[x_1, y]=0$ and $[x_1, a]=0$.
Since $a$ was arbitrarily chosen in $A$ we conclude that $[x_1, A]=0$. 
Let $D$ be the restricted ideal of $L$ generated by $A$ and $x_1$.  Then $D$ is abelian and yet $\dim L/D< \dim L/A$, which contradicts the choice of $A$. We deduce that $\dim A/N=1$, and then $N=Z(L)$. 

By virtue of the previous step, we can now replace $A$ with  $Z(L)+\F z$.

\item\label{[x1, x2]} 
\emph{ If $y_1, y_2 \in L$ are linearly independent modulo $Z(L)$ then  $[y_1, y_2]\neq 0$.}  In fact, if $[y_1, y_2]=0$,
then  $Z(L)+ \F y_1+\F y_2$ would be a restricted abelian ideal whose codimension in $L$ is less than the codimension of $A$ in $L$, a contradiction.

\item\label{L/Z(L) 3-dimensional}  \emph{$L/Z(L)$ is 3-dimensional.} If $\dim L/Z(L)\leq 2$, then $L$ contains an abelian restricted ideal of codimension at most 1, which is not possible. 
Suppose now that there exist $x_1, x_2, x_3, x_4\in L$ that are linearly independent modulo $Z(L)$.
Let  $z_{ij}=[x_i, x_j]$ for all $i,j$ and consider the element 
$$
v=[[x_4x_3x_1,x_4], [x_4x_1,x_1], x_2]=z_{14}^2\bigg( z_{12} z_{34} + z_{13}z_{24}+z_{14}z_{23}\bigg).
$$
By Theorem \ref{S-Z} and Lemma \ref{reduced}, we must have $v=0$.
Let $u= z_{12} z_{34} + z_{13}z_{24}+z_{14}z_{23}$. By symmetry we also get
$$
z_{24}^2u=z_{34}^2u=0.
$$
Thus
$$
z_{12}^2z_{34}^2u=z_{13}^2z_{24}^2u=z_{14}^2z_{23}^2u=0.
$$
It follows that
$$
\bigg(  z_{12} z_{34} + z_{13}z_{24}+z_{14}z_{23}\bigg)^3=0.
$$
Hence, by Lemma \ref{reduced} again, we deduce that 
\begin{align}\label{u=0}
 z_{12} z_{34} + z_{13}z_{24}+z_{14}z_{23}=0.
\end{align}
By Step \ref{[x1, x2]} we see that $z_{14}, z_{24}$ and $z_{34}$ are linearly independent. 
As a consequence, by Equation \eqref{u=0} and  the  PBW Theorem we deduce that 
each of $z_{12}$, $z_{13}$ and $z_{23}$ must be a linear combination of $z_{14}, z_{24},z_{34}$.
Since $ \la L^\prime\ra_2$ is free of $2$-nilpotent elements, we deduce that 
\begin{align*}
z_{12}&=\beta z_{14}+\gamma z_{24};\nonumber\\
z_{13}&=\alpha z_{14}+\gamma z_{34};\nonumber\\
z_{23}&=\alpha z_{24}+\beta z_{34}.\\
\end{align*}   
for some $\alpha, \beta, \gamma \in \F$.  Now put
 \begin{align*}
 y_1&=x_1+\gamma x_4;\\
y_2&=x_2+\beta x_4.
\end{align*}   
Then we have $[y_1, y_2]=0$, which contradicts Step \ref{[x1, x2]}. \qed 
\end{noindlist}

\begin{lem}\label{finite-3}
Let  $L$ be a   restricted Lie algebra over an algebraically closed field $\F$ of characteristic 2 such that $L^\prime$ is finite-dimensional and not 2-nilpotent. If $u(L)$ is Lie solvable then $L$ has  a finite-dimensional 2-nilpotent restricted ideal $I$ such that $\bar L=L/I$ satisfies one of the following conditions:
\begin{enumerate} 
\item[{\normalfont (i)}] $\bar L$ has  an abelian restricted ideal of codimension 1;
\item[{\normalfont (ii)}]  $\bar L$ is nilpotent of class 2 and $\dim \bar L/Z(\bar L)=3$;
\item[{\normalfont (iii)}]  $\bar L=  \la  x_1, x_2, y\ra_{\F} \oplus Z(\bar L)$, where  $[x_1,y]=x_1$,  $[x_2,y]=x_2$, and  $[x_1,x_2]\in Z(\bar L)$;
\item [{\normalfont (iv)}]  $\bar L=  \la  x, y\ra_{\F} \oplus H \oplus Z(\bar L)$, where  $H$ is a strongly abelian finite-dimensional restricted subalgebra of $\bar L$ such that $[x,y]=x$, $[y,h]=h$ and $[x,h]$ is a central element for every $h\in H$;
\item [{\normalfont (v)}] $\bar L= \la  x, y\ra_{\F} \oplus H \oplus Z(\bar L) $, where $H$ is a finite-dimensional abelian subalgebra of $\bar L$ such that $[x,y]=x$, $[y,h]=h$ and $[x,h]$ is a central element of $\bar L$ with $[x,h]^{[2]}=h^{[2]}$ for every $h\in H$.
\end{enumerate}
\end{lem}
\Proof 
First observe that, by arguing as in Step 1 in the proof of Lemma \ref{finite-2}, if $J$ is a finite-dimensional and 2-nilpotent restricted ideal of $L$ then we can replace $L$ with $L/J$. 
In particular, as the restricted ideal $\la[[L^\prime, L^\prime], L]\ra_2$ is finite-dimensional and 2-nilpotent by Theorem \ref{S-Z}, we can assume that $L^{\prime \prime}$ is central in $L$.
Moreover, by Theorem \ref{passman} there exists a 2-abelian restricted $A$ ideal  of $L$ of minimal finite codimention.   We replace $L$ with $L/ \la A^\prime \ra_2$ and thereby assume that $A$ is abelian. 
For every $x\in L$, let $J_x$ denote the restricted ideal generated by $[x,L]$.
Let $J=\sum J_x$, where  the sum runs over all $x\in L$ for which  $J_x\su Z(L)$
and $J_x$ is 2-nilpotent. Since $L^\prime$ is finite-dimensional, it is easy to see that $J$ is a 
finite-dimensional and 2-nilpotent restricted ideal of $L$. Thus  we can replace $L$ with $L/J$.
Let $I$ be the subspace consisting of all $z\in L^\prime\cap Z(L)$ such that $z$ is 2-nilpotent.
Since $\la I \ra_2$ is a finite-dimensional and 2-nilpotent restricted ideal of $L$, we can assume that $I=0$.
 Now, we consider two cases:

\begin{noindlistcase}
\item \emph{ Suppose $\zeta_2(L) \neq Z(L).$ }\label{Z2-neq-Z} Then there exists $z\in L$ such that $0\neq [z, L]\su Z(L)$.
 Note that $[z,L]$ is not 2-nilpotent. Thus, there exists $y\in L$ such that $[z, y]$ is not 2-nilpotent. Let $b\in A$, $x_1, x_2, x_3\in L$ and consider the following elements: 
 \begin{align}\label{uvw}\nonumber
u&=[[z,zy], [y, x_1y], x_2]=[z,y]^2[y, x_1, x_2];\\ \nonumber
v&=[[z,zy], [x_1, x_2y], b]=[z,y]( [x_1, x_2, b][z, y]+ [x_1, y, b][z, x_2]);\\ 
w&=[[z,zy], [x_1, x_2y], x_3]=[z,y]( [x_1, x_2, x_3][z, y]+ [y, x_1, x_3][z, x_2]).
\end{align} 
 Since, by Theorem \ref{S-Z},  $u$ is  nilpotent  and $[z,y]$ is not 2-nilpotent, we deduce that $[y, x_1, x_2]$ is 2-nilpotent.  
 Moreover, as $v$ and $[y,x_1, b]$ are both nilpotent we have that  $[x_1, x_2, b]$ is 2-nilpotent.
 As $x_1$ and $x_2$ were arbitrarily chosen, we conclude that $\la [L^\prime, A]\ra_2$ is a finite dimensional 2-nilpotent restricted ideal. 
Therefore we can replace $L$ by $L/\la [L^\prime, A]\ra_2$ and assume that $[L^\prime, A]=0$. We claim that $L^\prime \subseteq A$. Suppose the contrary. Regard $(L^\prime +A)/A$ as an $L$-module and let $\rho: L\to \text{End} ((L^\prime+A)/A)$ denote the corresponding representation. As $L^{\prime \prime} \subseteq Z(L)$ and $[L^\prime,A]=0$, the linear transformation $\rho([x, y])=\rho(x)\rho(y)-\rho(y)\rho(x)$ is nilpotent  for all $x, y\in L$. In view of  Theorem \ref{triang}, this entails that the linear Lie algebra $\rho(L)$ is triangularizable. Therefore $(L^\prime +A)/A$ contains a 1-dimensional $L$-module, that is to say there exists $z_1\in L^\prime$, $z_1\notin A$, such that $\la z_1\ra_2+A$ is a 2-abelian restricted ideal of $L$ of codimension less than $\dim L/A$, a contradiction.     
Thus, $L^\prime\su A$. In particular, $L^\prime$ is abelian and then, as $[x_1, x_2, x_3]$ and  $[y, x_1, x_3]$ commute and the elements $w$ and $[y, x_1, x_3]$ are both nilpotent, relation (\ref{uvw}) forces that  $[x_1, x_2, x_3]$ is nilpotent.
 As a consequence, the restricted ideal generated by $\gamma_3(L)$ is finite-dimensional and 2-nilpotent.
Hence, we can replace $L$ with $L/\la \gamma_3(L)\ra_2$ and assume that $L$ is nilpotent of class two. 
Then  Lemma \ref{finite-2} allows us to conclude that $L$ satisfies condition (i) or  (ii) of the statement.
 
\item \emph{Suppose $\zeta_2(L) = Z(L)$.} 
Clearly, $(L^\prime+Z(L))/Z(L)$ is a finite-dimensional $L$-module. Let 
$$
\rho: L\to \text{End} (L^\prime+Z(L)/Z(L))
$$
 denote the corresponding representation.  
Notice that, since $L^{\prime \prime}$ is central in $L$, the linear transformation $\rho([x, y])=\rho(x)\rho(y)-\rho(y)\rho(x)$ is nilpotent  for all $x, y\in L$. Therefore  Theorem \ref{triang} assures that the linear Lie algebra $\rho(L)$ is triangularizable.
Consequently, there exists a chain of ideals of $L$ 
$$
Z(L)=M_0\su M_1\su M_2\su \cdots \su M_n=L^\prime+Z(L)
$$
with  $\dim M_i/M_{i-1}=1$ for all $1\leq i\leq n$. We argue by induction on the triangularization length $n$ of the $L$-module  $(L^\prime+Z(L))/Z(L)$. If $n=0$ then $L^\prime\su Z(L)$ and so $L$ is nilpotent of class at most 2. The result then follows from Lemma \ref{finite-2}.  Now suppose that $n\geq 1$.
Since $\zeta_2(L) \subseteq Z(L)$,  $M_1/Z(L)$ is a 1-dimensional non-trivial $L$-module. Let $z_1\in M_1$, $z_1\notin Z(L)$. 
Then there exists  $y\in L$ such that  $[z_1, y]= z_1+ z_2$ for some $z_2\in Z(L)$. Put $z=z_1+z_2$, so that we have $[z, y]=z$. Let $\mathfrak{N}$ denote the annihilator of the $L$-module $M_1/Z(L)$. Clearly, $\fN$ is a restricted ideal of $L$ and $L={\F} y\oplus \fN$.
Now we consider two subcases:

\emph{Subcase 2.1: assume  $[z,\fN]=0.$}
Note that $\la z\ra_{2}$ is indeed a restricted  ideal of $ L$. If $z$ is  2-nilpotent then
we can replace $L$ with $L/\la z\ra_{2}$, so that the new $L$ has smaller triangularization length and the
result follows from the induction hypothesis. Hence, we assume that $z$ is not 2-nilpotent.
Let $b\in A$, $x\in \fN$, $x_1\in L$. Since $[z,\fN]=0$, it can be easily seen that $[z, b]=0$.  Consider 
 \begin{align*}
u&=[[y,z], [x_1, xy], b]=z[x_1, x, b].
\end{align*}
 Since $u$ is  nilpotent and $z$ is not 2-nilpotent, the element $[x_1, x, b]$ must be 2-nilpotent. 
Since $A$ is abelian and $L^\prime$ finite-dimensional, we  deduce that  $[L^\prime, A]$ is 2-nilpotent. Thus, we can replace $L$ with $L/\la [L^\prime, A] \ra_2$ and assume that $[L^\prime,A]=0$. 
In particular, $\la z\ra_2 + A$ is a 2-abelian restricted ideal of $L$ and so, by the minimality of the codimension of $A$, we must have $z\in A$.
Now, let $b\in A$, $x\in \fN$ and consider 
 \begin{align*}
v&=[[ by, z], [xb, x], y]=z [b, x]^2.
\end{align*}
 Since $v$ is  nilpotent and $z$ is not, the element $[b,x]$ must be  nilpotent. As a consequence, the restricted ideal $A+\fN$ is 2-abelian, which implies $\fN=A$.
 We conclude that $\dim L/A=1$, and condition (i) of the statement holds.
 
 \emph{Subcase 2.2: assume  $[z,\fN]\neq 0.$}  Let $x\in \fN$ such that $[x,z]\neq 0$.  Since $[x,z]\in L^\prime\cap Z(L)$, we deduce by the remarks prior to Case \ref{Z2-neq-Z} that  $[x,z]$ cannot be 2-nilpotent.
 Let $v_1,v_2\in \fN$ and $w\in L$. Similarly, if $[v_1,z]\neq 0$ then such an element is not 2-nilpotent. Consider the element 
 $$
 u=[[v_1,v_1v_2],[y,z],w]=[v_1,v_2,w][v_1,z].
 $$
 Since,  by Theorem \ref{S-Z}, $u$ is nilpotent and $[v_1,z]$ is not 2-nilpotent, the element $[v_1,v_2,w]$ is necessarily 2-nilpotent. On the other hand, 
 if  $[v_1,z]= 0$ then consider 
 $$
 v=[[xv_1,v_2],[y,z],w]=[v_1,v_2,w][x,z].
 $$
 As $v$ is nilpotent and $[x,z]$ is not 2-nilpotent, also in this case $[v_1,v_2,w]$ must be 2-nilpotent. As a consequence, the restricted ideal $\la [\fN^\prime,L]\ra_2$ is finite-dimensional and 2-nilpotent. Thus we can replace $L$ with $L/\la [\fN^\prime,L]\ra_2$ to assume that $\fN^\prime \subseteq Z(L)$.
 Let $f$ denote the vector space endomorphism induced by $\ad y$ on $\fN/Z(L)$. As $\zeta_2(L)=Z(L)$ and $f(\fN /Z(L))=([\fN,y]+Z(L))/Z(L)$ is finite-dimensional, we see that 
 $\ker f=0$ and $\dim \fN/Z(L)< \infty$.  
 Also, since $[z,y+y^{[2]}]=0$ we have $y+y^{[2]}\in \fN$. It follows that $f=f^2$ and then, as $\ker f=0$, $f$ acts the identity map on $\fN/Z(L)$.

We have $\fN=\fk\oplus Z(L)$, where $\fk$ denotes the subspace of $\fN$ consisting of the fixed points of $\ad y$. Of course, we can assume that $[\fk,a]\neq 0$  for every $a\in \fk$, $a\neq 0$, otherwise we can replace $z$ with $a$  and conclude by Subcase 2.1 that  condition (i) of the statement holds. 
Consider  
$$
J=\la x\in  [\fk, \fk]+\fk^{[2]} \mid x \mbox{ is 2-nilpotent} \ra_2.
$$
Since $\dim \fk< \infty$, $J$ is a central finite-dimensional 2-nilpotent restricted ideal of $L$. Hence we can replace $L$ with $L/J$.

\emph{Claim 1: For every $x_1,x_2\in \fk$ with $[x_1,x_2]\neq 0$ one has $C_{\fN}(x_1)\cap C_{\fN}(x_2)=Z(L)$.} Suppose that there exists   $0\neq b\in \fk$ such that $[x_1,b]=[x_2,b]=0$.  Since $[\fk, b]\neq 0$ there exists $c\in \fk$  such that $[b,c]\neq 0$. Consider the element
 $$
 \eta=[[x_1y,b], [y,x_2],c]=[x_1,x_2][b,c].
 $$ 
 Then $\eta$ is not nilpotent, which is not possible by Theorem \ref{S-Z}.  This proves the claim.
 
\emph{Claim 2: There exists no 3-dimensional subspace $P$ of $\fk$ with the property that for every basis $\{x_1,x_2,x_3\}$ of $P$ one has $[x_i,x_j]\neq 0$ whenever $i\neq j$.} Suppose otherwise and consider the element 
$$
\xi=[[y,x_1y],[y,x_2],x_3]=[x_1,x_2]x_3+[x_1,x_3]x_2+[x_2,x_3]x_1.
$$ 
Then $\xi$ is nilpotent by Theorem \ref{S-Z}. Note that by the remarks prior to Claim 1, the 
restricted Lie algebra $\la [x_i,x_j], x_k^2\mid 1\leq i, j, k\leq 3 \ra_2$ is free of nonzero 2-nilpotent elements. Thus, by Lemma \ref{reduced}, $\xi^2=0$.
Note that
\begin{align*}
\xi^{2}=[x_1,x_2]^{[2]}x_3^{[2]}+[x_1,x_3]^{[2]}x_2^{[2]}&+[x_2,x_3]^{[2]}x_1^{[2]}+[x_1,x_2][x_1,x_3][x_2,x_3].
\end{align*} 
However, the monomial $[x_1,x_2][x_1,x_3][x_2,x_3]$ has degree 3 and so, by the PBW Theorem, $\xi^{2}\neq 0$, a contradiction. This proves Claim 2.

\emph{Claim 3: $C_{\fN}(x)$ is abelian for every $x\in \fk$, $x\neq 0$.} Suppose otherwise. Then there exist $x_1,x_2\in C_{\fN}(x)\cap \fk$ such that $[x_1,x_2]\neq 0$. Thus, by Claim 1, we have $x\in C_{\fN}(x_1)\cap C_{\fN}(x_2)=Z(L)$, a contradiction. 

\emph{Claim 4: Let $x\in \fk$,  $x\neq 0$. If $\dim C_{\fN}(x)/Z(L)>1$ then $\dim {\fN}/ C_{\fN}(x)=1$.} There exists $a\in C_{\fN}(x)\cap \fk$ such that $a$ and $x$ are linearly independent modulo $Z(L)$. Let $b_1,b_2 \in \fk\setminus C_{\fN}(x)$. By Claim 1, we have $[a, b_1]\neq 0$ and $[a, b_2]\neq 0$.
  Now, Theorem \ref{S-Z} forces the nilpotency of the element
 $$
 \eta= [[xy,y], [a,y], b_1,b_2]=[x,b_1][a,b_2]+[x,b_2][a,b_1].
 $$ 
Lemma \ref{reduced} implies that $\eta=0$. Thus, by the PBW Theorem,
 $[x,b_1]$ must be proportional to at least one $[x,b_2]$ or $[a,b_1]$. But if $[x,b_1]=\beta [a,b_1]$ for some $\beta \in \F$, then $[x+\beta a,b_1]=0$. Thus,  $x+\beta a\in C_{\fN}(b_1)\cap C_{\fN}(x)$, which is not possible by Claim 1. As a consequence, we have $[x,b_1]=\alpha [x,b_2]$ for some $\alpha \in \F$. Hence,  
 $[x,b_1+\alpha b_2]=0$. We deduce that  $b_1+\alpha b_2\in C_{\fN}(x)$, and the claim follows.    
\emph{Claim 5: Let $x_1$ and $x_2$ be some   non-commuting elements of $\fk$. Then either 
$\dim C_{\fN}(x_1)/Z(L)=1$ or $\dim C_{\fN}(x_2)/Z(L)=1$.} Suppose that
 $\dim C_{\fN}(x_i)/Z(L)>1$, for $i=1, 2$. 
Since $[x_1, x_2]\neq 0$, we deduce by Claims 1 and 4 that 
$$1<\dim {C_{\fN}(x_1)}/Z(L)=\dim {C_{\fN}(x_1)}/({C_{\fN}(x_1)\cap C_{\fN}(x_2)}) \leq \dim {{\fN}}/{C_{\fN}(x_2)}=1,$$
a contradiction.

\emph{Claim 6: Suppose that there exists $0\neq x_1\in \fk$ such that $\dim C_{\fN}(x_1)/Z(L)>1$. Then either condition (iv) or condition (v) of the statement holds.} 
Note that, by Claim 5, $\dim C_{\fN}(x)/Z(L)=1$, for every  $x\in \fk\backslash {C_{\fN}(x_1)}$.
Recall that if $b\in \fk$ is 2-nilpotent then $b^{[2]}=0$. 
By Claim 3, we know that $C_{\fN}(x_1)$ is  abelian. Also,
if $x_2, x_3\in {\fN}$ are linearly independent modulo $C_{\fN}(x_1)$ then, by Claim 2, $[x_2, x_3]=0$.
But then $\dim C_{\fN}(x_2)/Z(L)\geq 2$, which is a contradiction.
Thus, $\dim {\fN}/C_{\fN}(x_1)=1$. Now we take $H= \fk \cap {C_{\fN}(x_1)}$.
Let $x\in \fk\backslash H$ and  $h\in H$ such that $h$ and  $x_1$ are not proportional.
 Now consider
$$
\zeta=[[x_1y,y],[h,y],x]=[x,x_1]h+[x,h]x_1.
$$    
Then $\zeta$ is nilpotent by Theorem \ref{S-Z}. Now we apply Lemma \ref{reduced} to the restricted Lie algebra $\la [x,x_1], h^{[2]}, [x,h], x_1^{[2]}\ra_{2}$. We deduce that $\zeta^2=0$.
Note that  $[x,x_1]^{[2]}\neq 0$ and $[x,h]^{[2]}\neq 0$ because $[x,x_1]$ and $[x,h]$ are central and nonzero. Consequently, by the PBW Theorem, if $x_1^{[2]}=0$ then  $h^{[2]}=0$.
We conclude that if $x_1^{[2]}=0$  then $H$ is strongly abelian,  
and  condition (iv) of the statement holds. Now suppose that $x_1^{[2]}\neq 0$. Observe that if $h^{[2]}=\alpha x_1^{[2]}$ for some $\alpha \in F$ then the element $h+\alpha^{\frac{1}{2}} x_1$ is 2-nilpotent, thus by replacing $h$ with $h+\alpha^{\frac{1}{2}}x_1$ we see that $\zeta$ is not nilpotent, a contradiction to Theorem \ref{S-Z}.  Hence $x_1^{[2]}$ and $h^{[2]}$ are linearly independent. Since $\zeta^2=0$, for some $\beta \in {\F}$,  $\beta \neq 0$, we must have
$$
[x,x_1]^{[2]}=\beta x_1^{[2]}, \quad [x,h]^{[2]}=\beta h^{[2]}.
$$ 
Put $\alpha=\beta^{\frac{1}{2}}$.
We replace $x$ with $\alpha x$ to assume that $[x,h]^{[2]}=h^{[2]}$, and condition (v) of the statement holds.

\emph{Claim 7: If $\dim C_{\fN}(x)/Z(L)=1$, for every $x\in \fk$, $x\in \fk$, then condition (iii) of the statement holds.} By Claim 2, it is clear that in this case $\dim {\fN}/Z(L)=2$, thereby condition (iii) of the statement holds. This completes the proof.
  \qed
\end{noindlistcase}

\section{Infinite-dimensional derived subalgebra}

In this section we handle restricted Lie algebras with derived subalgebra of infinite dimension.

\begin{lem}\label{infinite-xy}
Let $L$ be a restricted Lie algebra over a field $\F$ of  characteristic 2. Suppose that $L$ contains an abelian restricted ideal $I$ of codimension 2 such that $L/I$ is $2$-nilpotent. If $[x, I]$ is infinite-dimensional for every  $x\in L\backslash I$ then $u(L)$ is not Lie solvable.
\end{lem}
\Proof Suppose  the contrary. There exist $x,y\in L$ that are linearly independent modulo $I$ and such that $y^{[2]}\in I$.  Clearly,  $L/I$ is abelian, in particular $[x,y]\in I$. Moreover, for all $a\in I$ we have 
$[a, x, y]=[a, y, x]$ and so $[[I,x],y]=[[I,y],x]$. Since $L/I$ is 2-nilpotent,  $x^{[2]}=\beta y$ modulo $I$, for some $\beta\in \F$. It follows that $[[I, x], y]\su Z(L)$.
Now we consider the following cases:

\begin{noindlistcase}

\item \emph{ $[[I,x],y]$ is finite-dimensional.} Then there exists a subspace $A$ of finite codimension in $I$ such that $[[A,x],y]=0$. We can replace $I$ with $A$ to assume that 
 $[[I,x],y]=0$. Let us consider two subcases. 

 \emph{Subcase 1.1: There exists a sequence $b_1,b_2,\ldots$ of elements of $I$ such that the $[b_i,y]$ are linearly independent and the subspace $V$ spanned by the $[b_i,x]$ has finite dimension $t$.} Clearly, we can assume that $[b_1, x],\ldots, [b_t, x]$ span $V$.
We can then rescale the $b_{i}$ in such a way that   $[b_i,x]=0$, for all $i\geq t+1$. We discard $b_1, \ldots, b_t$ and relabel $b_{t+1}, b_{t+1}, \ldots$ to $b_1, b_2, \ldots$.
Since $[I,x]$ is infinite-dimensional,  for every even integer $N=2k$ there exist  $a_1,a_2,\ldots, a_k \in I$  such that the set consisting of all of the  $[x, a_i]$ and $[y, b_j]$, $i\leq k, j\leq N$, is linearly independent.
For every even integer $i$  consider the elements
\begin{align*}
A_i&=[xa_ib_i,x]=x[a_i,x]b_i;\\
B_i&=[yb_{i+1},y]=y[b_{i+1},y].
\end{align*}
Then we have
$$
[A_i,B_i]=x[a_i,x][b_i,y][b_{i+1},y] + [x,y][a_i,x][b_{i+1},y]b_i
$$
and so
\begin{align*}
C_i=[[A_i,B_i],a_{i+1}]=[a_i,x][a_{i+1},x][b_i,y][b_{i+1},y].
\end{align*}
By the PBW Theorem,    for every even integer $N$ there exists a nonzero
element $C_2C_4\cdots C_{N}\in [[u(L), u(L)], [u(L), u(L)], u(L)]^k$.
Therefore we conclude by Theorem \ref{S-Z} that $u(L)$ is not Lie solvable.

Note that the just considered subcase occurs in particular when $y=\beta x^{[2]}$ modulo $I$,  for some $\beta \in \F$, $\beta \neq 0$. Moreover, as $x$ and $y$ are linearly independent modulo $I$ and $y^{[2]}\in I$, it is not possible to have $y=\alpha x+ \beta x^{[2]}$ modulo $I$ for some $\alpha, \beta \in \F$. Therefore   in the next subcase we can assume that $x^{[2]}\in I$, as well.  

\emph{Subcase 1.2: For every sequence $b_1,b_2,\ldots$ of elements of $I$ such that the $[b_i,y]$ are linearly independent one has that  the subspace $V$ spanned by the $[b_i,x]$ is infinite dimensional.}   
We consider two cases.

\emph{ 1.2.1:} Suppose that for every integer $k$ there exists   $a_1, \ldots, a_k\in A$ 
such that the subspace spanned by all of the $[x, a_i]$ and  $[y, a_j]$, $1\leq i, j\leq k$, has dimension $2k$. Let $k$ be an even integer. 
Let $i$ be an  odd integer in the range $1\leq i \leq k$ and  consider 
\begin{align*}
C_i&=[[xa_{i},ya_{i+1}], [y,xa_{i}], a_{i+1}]\\
&=[x,a_{i}] [y, a_{i}] [x,a_{i+1}][y,a_{i+1}] +[x, a_{i+1}]^2[y,a_{i}]^2.
\end{align*}
Notice that all monomials in the product $C_1C_3\cdots C_{k-1}$ have degree less than $2k$ except for  
$[x,a_1]\cdots[x, a_k] [y, a_1] \cdots [y,a_k]$. Consequently, by the PBW Theorem we have $C_1C_3\cdots C_{k-1}\neq 0$  and then, by Theorem \ref{S-Z}, $u(L)$ cannot be Lie solvable.

\emph{1.2.2: Suppose that  1.2.1 fails.} 
 Let $n$ be the largest integer such that  1.2.1 holds.
Then there exists $a_1, \ldots, a_n\in I$  such that the subspace $S$  spanned by all of the  $[y, a_i]$ and $[x, a_j]$, $1\leq i, j\leq n$, has dimension $2n$. Since  1.2.1 fails for $n+1$,  for every    $b\in I \backslash \langle a_1,\ldots,a_n\rangle_{\F}$ there exists $\alpha_b\in \F$ such  that    $[x, b] =\alpha_b [y, b]$ modulo $S$. 
Let $D$ be the set consisting of all $\alpha_b$ and for every $\alpha \in D$ denote by $I_{\alpha}$ the subspace of $A$ consisting of all $b$ such that $[x, b] =\alpha [y, b]$ modulo $S$. Then for every $\alpha \in \F $ we have that $I/I_\alpha$ is infinite dimensional. In fact, in the contrary case we would have $\dim [x+\alpha y,I]<\infty$, a contradiction. 
As $S$ is finite-dimensional, if the set $D$ is finite then there exist distinct nonzero
$\alpha_1, \alpha_2\in \F$ such that  $I_{\alpha_1}$ and $I_{\alpha_2}$ are both infinite-dimensional.

Now let  $b_1, b_2, \ldots$ in $I$ such that $[x, b_1], [x, b_2], \ldots $ are linearly independent. For every positive integer $i$ write $\alpha_i$ for 
$\alpha_{b_i}$.
 If $D$ is infinite then without loss of generality we can assume that $\alpha_i\neq \alpha_j $ whenever $ i\neq j$.  
 On the other hand,  if $D$ is finite then for every $t>0$ we  choose $b_1, \ldots, b_{4t}$ so that $\alpha_i\neq \alpha_{i+3}$, for all $1\leq i\leq 4t-3$. 
For every $0\leq i \leq 4t$ put $u_i=[xb_{i},yb_{i+1}]$ and $v_i=[y,xb_{i+2}]$. 
We have
$$
C_i=[[u_i,v_i],b_{i+3}]=\alpha_{i+2} (\alpha_{i+3}+\alpha_{i}) [x,b_{i}] [x, b_{i+1}]  [x,b_{i+2}] [x, b_{i+3}] +w_i,
$$
where $w_i$ is a linear combination of PBW monomials that each involve at least an element of $S$. Now consider 
 the product $C=C_1\cdots C_{4i}\cdots C_{4t}$. We observe that $C=\alpha [x,b_{1}] \cdots [x, b_{4t}]+w$, where $0\neq\alpha\in \F$ and each PBW monomial involved in $w$ either  has degree less than $4t$ or involves at least an element of $S$. From the PBW Theorem it  follows that $C\neq 0$, and by Theorem \ref{S-Z} we conclude that $u(L)$ is not Lie solvable.
 
\item \emph{$ [[I,x],y]$ is infinite dimensional.} Notice that $[[I,x],y]$ is contained in $Z(L)$. For every $n$ positive integer $n$ pick $a_1,a_2,\ldots, a_n \in [I,y]$ and $b_1,b_2,\ldots, b_n \in [I,x]$ such that the set $\{[a_i,x], [b_j,y]\, \vert \, 1\leq i,j\leq n \}$ is linearly independent. Since $y$ and the $a_i$ commute, we have
$$
D_i=[[xa_{2i-1}b_{2i-1},a_{2i}],[y,yb_{2i}], x]=[x,a_{2i-1}][x,b_{2i-1}][x,a_{2i}][x,b_{2i}].
$$
By the PBW Theorem we have $D_1D_2\cdots D_n \neq 0$, thus Theorem \ref{S-Z} forces that $u(L)$ is not Lie solvable.   
This finishes the proof. \qed
\end{noindlistcase}

\begin{lem}\label{infinite-xynonnilp}
Let $L$ be a restricted Lie algebra over an algebraically closed field $\F$ of characteristic 2. Let $I$ be an abelian restricted ideal of $L$ of codimension 2 such that $L^\prime \su I$ and there exists a nonzero toral element in $L/I$. If $[x, I]$ is infinite-dimensional for every  $x\in L\backslash I$ then $u(L)$ is not Lie solvable.
\end{lem}
\Proof
By assumption there exists $y\in L\backslash I$ such that $y^2=y$ modulo $I$. Let $x\in L$ such that $x$ and $y$ are linearly independent modulo  $I$. Then, by hypothesis, $ [I,x] $ and $[I,y]$ are both infinite dimensional. Moreover, by Jacobi identity we have $[[I,x],y]=[[I,y],x]$. Let us consider two cases separately. 

{\bf Case 1.} \emph{$\dim [[I,x],y]<\infty$.}  Put $B=[I,y]$. As $\dim [[I,x],y]<\infty$, we can find  a sequence $a_1,a_2,\ldots$ of elements of $I$ such that the $[a_i,x]$ are linearly independent and $[a_i,x,y]=0$ for every $i$. Moreover, as $[B,x]=[[I,y],x]$ is finite-dimensional,   there exist linearly independent elements $b_1, b_2, \ldots $  in  $B$ such that $[b_i,x]=0$, for all $i$.
Note that, since $y^{[2]}=y$ modulo $I$, we have $[b,y^{[2]}]=[b,y]$ for every $b\in I$. 
Thus,   the set $\{[a_i,x],b_j \vert \, i,j>0 \}$ is linearly independent. For every even integer $i$,  consider the element
$$
C_i=[[xa_ib_i,x],[yb_{i+1},y],a_{i+1}]=[a_i,x][a_{i+1},x]b_ib_{i+1}.
$$
By the PBW Theorem,    for every even integer $N$ there exists a nonzero
element $C_2C_4\cdots C_{N}\in [[u(L), u(L)], [u(L), u(L)], u(L)]^N$.
Therefore, by Theorem  \ref{S-Z}  we conclude that $u(L)$ is not Lie solvable. 

{\bf Case 2.} \emph{$\dim [[I,x],y]=\infty$.} We split this case in two subcases.

\emph{Subcase 2.1: The power mapping of $L/I$ is singular.}  
We may assume  that of $x^{[2]}\in I$.
If $\dim [[I,x],y]<\infty$ then $u(L)$ is not Lie solvable by Case 1. Thus we assume that $\dim [[I,x],y]=\infty$. 
There exist $a_1,a_2,\ldots$ in $I$  such that  the $[a_i,x,y]$ are linearly independent.
Let $b_i=[a_i, x,y]$. Note that $[b_i, y]=b_i$, $[a_i, y^{[2]}]=[a_i, y]$, and $[b_i, x]=0$, for all $i$.
For every $k$ consider the element
$$
B_k=[[a_{k}, y], [y,yb_{k+1}],x]=[[a_k, y]b_{k+1}, x]=b_kb_{k+1}.
$$
Then the PBW Theorem yields that $B_1B_3\cdots B_{2k+1}\neq 0$  for every $k>0$, contradicting Theorem \ref{S-Z}.

\emph{Subcase 2.2: The power mapping of $L/I$ is nonsingular.} As the ground field is algebraically closed, by \cite[Chapter 2, \S 3, Theorem 3.6]{SF} the restricted Lie algebra $L/I$ has a toral basis, in particular we can assume that $x^{[2]}=x$ modulo $I$.  
Since $\dim [x+y,I]=\infty$,  there exist  $a_1,a_2,\ldots$ in $I$ such that the $[a_i,x+y]$ are linearly independent.  Also, since $\dim [I, x, y]=\infty$,  there exist  $b_1,b_2,\ldots$ in $I$
such that the $c_j=[b_j, x, y]$ are linearly independent. Note that $[c_j, x]=[c_j, y]=c_j$, for all $j$.
Thus, $[c_j, x+y]=0$. We deduce that the set $\{[a_i,x+y], c_j\, \vert \, i,j>0 \}$ is linearly independent.  For every $i>0$ consider the element
$$
C_i=[[[x,xc_{2i-1}],[y,yc_{2i}]],a_i]=[[xc_{2i-1},yc_{2i}],a_i]=c_{2i}c_{2i-1}[a_i,x+y].
$$
By the PBW Theorem we have that $C_1C_2\cdots C_n\neq 0$ for every $n>0$, therefore $u(L)$ is not Lie solvable by  Theorem \ref{S-Z}, yielding the claim. \qed

\begin{lem}\label{infinite-metab}
Let $L$ be a restricted Lie algebra over an algebraically closed field $\F$ of characteristic 2.
Suppose that $u(L)$ is Lie solvable and let $I$ be a 2-abelian restricted ideal of $L$ of minimal  codimension. If $L^\prime$ is infinite dimensional and $L^\prime\su I$ then $\dim L/I \leq 1$.
\end{lem}

\Proof  Suppose, if possible, that $\dim L/I\geq 2$. We replace $L$ with $L/\la I^\prime \ra_2$ to assume that $I$ is abelian.
\begin{noindlistcase}
\item
\emph{There exists $x\in L\backslash I$ such that $\dim[I,x]$ is finite.}
By the maximality of $I$ we see that $[x,a]$ is not 2-nilpotent, for some $a\in I$. 
 As $L^\prime$ is infinite dimensional, there exists $y \in L$ such that $ [y,I]$ is infinite-dimensional. 
 By Theorem \ref{S-Z} there exists $N>0$ such that  $[[[u(L),u(L)],[u(L),u(L)]],u(L)]^N=0$.  Let $k$ be an integer such that $2^{k-1}\geq N$ and pick $b_1, b_2, \ldots,b_{2^{k}} \in I$ such that the set of all $[y, b_i]$ and the element $[x,a]^{[2]^{k}}$ are linearly independent. As $[x, I]$ is finite dimensional  the $b_i$ can be chosen so that $[x, b_i]=0$, for all $i$.
For every $i=1,2,\ldots,2^{k-1}$ put
$$
u_{2i-1}=[ [axy, b_{2i-1}], [b_{2i},x y],  a ]=[x, a]^2[y, b_{2i-1}][y, b_{2i}].
$$
Then we have 
$[x, a]^{[2]^{k}}[y, b_1]\cdots [y, b_{2^k}]=u_1u_3\cdots u_{2^k-1}=0$, which  contradicts  the PBW Theorem. 

\item \emph{For every  $x\in L\backslash I$ one has that $[x, I]$ is infinite-dimensional.}
   As $L/I$ is abelian and finite-dimensional (by Theorem \ref{passman}) and, moreover, the ground field is algebraically closed, by \cite[Chapter 4, Theorem 4.5.8]{W} we have $L/I=\mathfrak{T}\oplus \mathfrak{N}$, where $\mathfrak{T}$ is a torus and $\mathfrak{N}$ is a 2-nilpotent restricted subalgebra. In turn, by \cite[Chapter 2, \S 3, Theorem 3.6]{SF},  $\mathfrak{T}$ has a toral basis. In particular, $L/I$ contains a 2-dimensional subalgebra  which is 2-nilpotent or contains a nonzero toral element. 
 The result then follows from  Lemmas  \ref{infinite-xy} and  \ref{infinite-xynonnilp}. 
\qed
\end{noindlistcase}

\begin{lem}\label{infinite-4}
Let $L$ be a restricted Lie algebra over an algebraically closed field $\F$ of characteristic 2 such that $L^\prime$ is infinite-dimensional. If  $u(L)$ is Lie solvable then $L$ has a 2-abelian restricted ideal of codimension 1.
\end{lem}
\Proof 
Note that if $J$ is a finite-dimensional and 2-nilpotent restricted ideal of $L$, then we can replace $L$ with $\bar L=L/J$. 
In particular, as $\la [[L^\prime, L^\prime], L]\ra_2$ is finite-dimensional and 2-nilpotent by Theorem \ref{S-Z},  we can assume that $L^{\prime \prime}$ is central in $L$.
Moreover, by Theorem \ref{passman} there exists a 2-abelian restricted ideal  of $L$ of finite codimention. Let  $I$ be a such restricted ideal of $L$ of minimal codimension.  
In view of Lemma \ref{infinite-metab}, in order to prove the statement it is enough to show that $L^\prime\su I$. Suppose the contrary. We can replace $L$ with $L/\la I^\prime \ra_2$ and thereby assume that $I$ is abelian. 
 
Consider $L^\prime +I/I$ as an $L$ module and let $\rho: L\to \text{End} (L^\prime +I/I)$ denote the corresponding representation.  
As $L^{\prime \prime}$ is central in $L$, $\rho([x, y])=\rho(x)\rho(y)-\rho(y)\rho(x)$ is a nilpotent transformation on $L^\prime$ for all $x, y\in L$. Thus, by Theorem \ref{triang}, the linear Lie algebra $\rho(L)$ is triangularizable. 
Consequently, there exists a chain of ideals of $L$ 
$$
I= M_0\su  M_1\su  M_2\su \cdots \su  M_n= L^\prime + I
$$
such that $\dim  M_i/M_{i-1}=1$ for all $i=1,2,\ldots,n$.
In particular,  there exists a 1-dimensional ideal of $L/I$ of the form $M=(\F z+I)/I$ for a suitable $z\in L^\prime$. Also, by Theorem \ref{S-Z}, there exists a positive integer $n$ such that 
$$[ [u(L),u(L)],[u(L),u(L)],u(L)]^n=0.$$ 
Now we consider two cases:
\begin{noindlistcase}
\item 
 \emph{Suppose  $[L,z]\nsubseteq I$.} Then there exists  $y\in L$ such that  $[z, y]= z$ modulo $I$. 
 It follows that   $z^{[2]}\in I$. 
 Let $\fN$ denote the annihilator of the $L$-module $M$. Note that 
 $\fN$ is a restricted ideal of $L$ and $L=\F y\oplus \fN$. Suppose, if possible, that $\dim [I,z]<\infty$. Since $C=\la z\ra_2 +I$ is a restricted ideal of $L$ and $\dim L/C<\dim L/I$, by the minimality of $\dim L/I$ there exists $b\in I$ such $[z,b]$ is not 2-nilpotent. As $L^\prime$ is infinite dimensional, there exists $x\in L$ such that $\dim [I,x]=\infty$. 
 Let $b_1,\ldots,b_{n}\in I$ such that the set of all $[x,b_i]$ and $[z,b]^{[2]^n}$ is linearly independent.  As $\dim [I,z]<\infty$, the $b_i$  can be chosen so that $[z,b_i]=0$, for every $i$. Put 
$$
v_i=[[z,y], [b_i, zbx], b]=[z, b]^2[x,b_i].
$$
Then  $v_1v_2\cdots v_{n}=[z, b]^{2^n}[x,b_1]\cdots[x,b_{n}]=0$, contradicting the PBW Theorem. Thus 
$[I,z]$ is infinite-dimensional.  Note that $[I,z]=[I,[y,z]]\subseteq [[I,z],y]+[[I,y],z]]$, in particular  $[[I,y],z]]$ and $[[I,y],z]]$ cannot both be finite-dimensional. Let us split this case in some subcases.  

\emph{Subcase 1.1: $\dim[[I,y],z]]<\infty$.} Then there exist $a_1,  \ldots, a_n \in I$ such that the $[a_1,y],\ldots, [a_n,y]$ are linearly independent and $[[a_i,y],z]=0$,  for every $1\leq i\leq n$. Now take $b_1,\ldots, b_n\in I$ such that the set of all $[a_i,y]$ and $[b_i,z]$ is linearly independent. For every $i$, consider the element
$$
u_i=[[y,ya_i],[y,z], b_i]=[a_i,y][b_i,z].
$$ 
Therefore we must have  $u_1u_2\cdots u_n=0$. But this contradicts the PBW Theorem.

\emph{Subcase 1.2: $\dim[[I,y],z]]=\infty$ and $\dim[[I,z],y]]<\infty$.} 
There exist $a_1, \ldots, a_{2n}\in I$ such that $[a_1,z], \ldots, [a_n,z]$ are linearly independent and $[a_i,z,y]=0$,  for every $1\leq i\leq n$. Note that $[a_i,y,z]=[a_i, z]$,  for every $1\leq i\leq n$. 
Since $z^{[2]}\in I$, we have 
$$
u_i=[[z, y],[a_i,a_{i+1}y], z]= [z, a_{i+1}][a_i, y, z]=[z, a_{i+1}][z, a_i].
$$ 
Then we must have  $u_2u_4\cdots u_{2n}=0$, contradicting the PBW Theorem.

\emph{Subcase 1.3: $\dim[[I,y],z]]=\infty$ and $\dim[[I,z],y]]=\infty$.}
Then there exist $a_1, \ldots, a_n\in I$ and $b_1,\ldots, b_n\in I$  such that the set of all 
$[a_i,z, y]$ and $[b_j,y, z]$   is linearly independent. 
Since $[b_i, y, z]$ is a central element, we have  
$$
u_i=[[a_iz,z],[b_i,y],y]=[a_i,z,y][b_i,y,z].
$$ 
Thus $u_1u_2\cdots u_n=0$,  which contradicts the PBW Theorem, again.

 \item 
 \emph{Suppose  $[L,z]\su I$.} We distinguish two subcases:
 
\emph{Subcase 2.1: $\dim[I,z]<\infty$.}  Since  $\dim L/I$ is minimal, there exists $a\in I$ such that $[z,a]$ in not 2-nilpotent. Moreover, as $L^\prime$ is infinite dimensional, there exists   $x \in L$ such that $\dim [I,x]=\infty$.
Let  $b_1,  \ldots,b_{n}\in I$ such that the set of all $[x, b_i]$ and the element  $ [z, a]^{[2]^{n}}$ is linearly independent. As $[I,z]$ is finite-dimensional, we can choose the $b_i$  so that $[b_i, z]=0$, for all $i$. 
Since $[z, x]\in I$, we also have $[b_i, x, z]=0$, for all $i$. For every $1\leq i\leq n$, consider   
$$
u_{i}=[[ x, yb_{i}], z,  a]=[z, a][x, b_{i}].
$$
Then we have  
$u_1u_2\cdots u_{n}=[z, a_1]^{{n}}[x, b_1]\cdots [x, b_{n}]=0$, which is a contradiction to the PBW Theorem. 

\emph{Subcase 2.2: $\dim[I,z]=\infty$.} Let $x_1,y_1,\ldots,x_n,y_n\in L$ and nonzero $\alpha_1,\ldots,\alpha_n\in \F$ such that $z=\sum_{i=1}^n\alpha_i[x_i,y_i]$. From the assumption it follows that for some $i$ the subspace $[I,[x_i,y_i]]$ is infinite dimensional. 
 Let $H$ be the restricted subalgebra generated by $I$, $x_i$ and $z$. Note that for every $\alpha \in \F$ one has 
 $$
 [I,[\alpha x_i,y_i]]=[I,[\alpha x_i+\beta z,y_i]]\subseteq [[I,\alpha x_i+\beta z],y_i]]+[[I,y_i],\alpha x_i+\beta z]].
 $$
 Note that if $\dim [I,\alpha x_i+\beta z]$ is finite then so is 
  $\dim [[I,y_i],\alpha x_i+\beta z]]$, which would imply that $\dim [I,[\alpha x_i,y_i]]$ is finite.
We deduce that $\dim [I,\alpha x_i+\beta z]=\infty$ unless $\alpha=0$.
Consequently, for every $\alpha,\beta \in \F$, $(\alpha,\beta)\neq (0,0)$, one has $\dim[\alpha x_i+\beta z,I]=\infty$. It follows that $H$ cannot contain any 2-abelian restricted ideal of codimension 1. Therefore, by Lemma \ref{infinite-metab}, we conclude that $u(H)$ (and then $u(L)$) is not Lie solvable, a contradiction, and the proof is complete. \qed
\end{noindlistcase}

\section{Proof of the Main Theorem} 
In the next lemmata we show that each of the conditions (i)-(v) in the statement of the Main Theorem implies Lie solvability of $u(L)$. 

\begin{lem}\label{cond-i}
Let $L$ be a restricted Lie algebra over a field of characteristic 2. If $L$ contains a 2-abelian restricted ideal $I$ of codimension at most 1 then  $u(L)$ is Lie solvable.
\end{lem}
\Proof 
If $I=L$ then $u(L)$ is Lie solvable by \cite{RS1}. Suppose then that $I$ has codimension 1 in $L$ and put $A=I/\la I^\prime\ra_2$. Then $A$ is an abelian restricted ideal of $L/\la I^\prime \ra_2$ and $u(L/\la I^\prime \ra_2)$ is a free left $u(A)$-module of rank 2. As a consequence, it follows easily that $u(L/\la I^\prime \ra_2)$ embeds in a matrix algebra $M_2(u(A))$. Moreover, as $u(A)$ is commutative and $\F$ has characteristic 2, the algebra $M_2(u(A))$ is Lie solvable. Thus  $u(L/\la I^\prime \ra_2)\cong u(L)/u(L)\la I^\prime \ra_2$ is Lie solvable and so, as the ideal $u(L)\la I^\prime \ra_2$ is nilpotent, we conclude that $u(L)$ is Lie solvable.
\qed

\begin{lem}\label{cond-ii}
Let $L$ be a  nilpotent restricted Lie algebra of class 2 over a field of characteristic 2. If $\dim L/Z(L)\leq 3$ then $u(L)$ is Lie solvable.
\end{lem}
\Proof Let $x_1, x_2, x_3\in L$ be linearly independent modulo $Z(L)$.
Consider  $\mathfrak{g}=u(L)$ as a Lie algebra and put $\mathfrak{Z}=u(Z(L))$. We will prove that $\mathfrak{g}$ is solvable.
For this purpose, consider 
$$
\mathfrak{H}=\la z_1x_1, z_2x_2, z_3x_3, z_4x_1x_2, z_5x_1x_3, z_6x_2x_3, z_7\mid z_1, \ldots, z_7\in \mathfrak{Z}\ra_{\F}.
$$
Observe  that $[\mathfrak{H}, \mathfrak{g}]\su \mathfrak{H}$. Since $\fg/\fh$ is abelian, it is enough to prove that 
$\fh$ is solvable. Now let
$$
\mathfrak{k}=\la z_1x_1, z_2x_2, z_3x_3, z_4 \mid z_1, \ldots, z_4\in \mathfrak{Z}\ra_{\F}.
$$
Note that $[\fk, \ \fh]\su \fk$. Since $\fk^{\prime \prime}=0$, it is enough to prove that $\mathfrak{n}=\fh/\fk$ is solvable.
We have  $\mathfrak{n}=(\fm + \fk)/\fk$, where $\fm=\la z_1x_1x_2, z_2x_1x_3, z_3x_2x_3\mid z_1, z_2, z_3\in \mathfrak{Z}\ra_{\F}$.
 Put $e_1=[x_1x_2, x_1x_3], e_2=[x_1x_2, x_2x_3], e_3=[x_1x_3, x_2x_3]$.
We have
\begin{align*}
e_1&=[x_1, x_3]x_1x_2+ [x_1, x_2] x_1x_3\qquad \textrm{modulo} \fk,\\
e_2&=[x_1, x_2]x_2x_3+ [x_2, x_3] x_1x_2 \qquad \textrm{modulo} \fk,\\
e_3&=[x_1, x_3]x_2x_3+ [x_2, x_3] x_1x_3 \qquad \textrm{modulo} \fk.
\end{align*}
At this stage, a simple calculation shows that $[e_i, e_j]\in \fk$, for all $1\leq i, j\leq 3$. We conclude that $\mathfrak{n}^{\prime \prime}=0$, yielding the claim.
\qed

\begin{lem}\label{cond-iii} Let $L= \la  x_1, x_2, y\ra_{\F} \oplus Z( L)$ be a restricted Lie algebra over a field $\F$ of characteristic 2, where  $[x_1,y]=x_1$,  $[x_2,y]=x_2$,  $[x_1,x_2]\in Z( L)$, and the power mapping is arbitrary. Then $u(L)$ is Lie solvable.
\end{lem}
\Proof We can suppose that $x_1, x_2, y$ are linearly independent (otherwise $L$ contains an abelian restricted ideal of codimension at most 1 and the claim follows from Lemma \ref{cond-i}).
Consider  $\mathfrak{g}=u(L)$ as a Lie algebra and put $\mathfrak{Z}=u(Z(L))$. We will  prove that $\mathfrak{g}$ is solvable.
For this purpose,  consider 
$$
\mathfrak{H}=\la z_1x_1, z_2x_2, z_3x_1y, z_4x_2y, z_5\mid z_1, \ldots, z_5\in \mathfrak{Z}\ra_{\F}.
$$
 Note that $x_1^2, x_2^2, [x_1, x_2y]\in \mathfrak{Z}$, $ [x_1y, x_2y]=0$, and $y^{[2]}=y$ modulo $Z(L)$. It follows that   $[\mathfrak{H}, \mathfrak{g}]\su \mathfrak{H}$. Since $\fg/\fh$ is abelian, it is enough to show that
 $\fh$ is solvable. One has that
 \begin{align*}
 \mathfrak{H}^\prime\su \la z_1x_1, z_2x_2, z_3\mid z_1, \ldots, z_3\in \mathfrak{Z}\ra_{\F}.
 \end{align*}
 It is easy now to see that $\mathfrak{H}^{\prime \prime}\su \mathfrak{Z}$. Hence, $\mathfrak{H}$ is solvable, as required.
\qed

\begin{lem}\label{cond-iv} Let $L = \la  x , y\ra_{\F}  \oplus H \oplus Z(L)$ be a restricted Lie algebra over a field $\F$ of characteristic 2 that satisfies either  condition (iv) or (v) of the Main Theorem. Then $u(L)$ is Lie solvable.
\end{lem}
\Proof
We can suppose that $x$ and $y$ are linearly independent (otherwise $u(L)$ is Lie solvable by Lemma \ref{cond-i}). 
Note that $x^{[2]}\in H+Z(L)$ and $y^{[2]}=y$ modulo $Z(L)$. Let us evaluate the commutators.
Let $h_1, \ldots, h_n$ be a basis of $H$. Set $\fz=u(Z(L)$ and $R=u(H+Z(L))$.
Note that elements of $R$ are linear combinations of PBW monomials of the form 
$zh_1^{\epsilon_1}\ldots h_n^{\epsilon_n}$, where $z\in \fz$ and $\epsilon_i\in \{0, 1\}$.
Let $\bZ_0$ and $\bZ_1$ be the set of  even and  odd integers respectively.
We define the degree $\deg_H$ of any such monomial as $\epsilon_1+\cdots+\epsilon_n$.
Let $u$ be a PBW monomial in $R$. Note that if $\deg_H (u)\in \bZ_0$  then $[u, y]=0$.
Also, if $\deg_H (u)\in \bZ_1$ then $[u, y]=u$. 
Now, let $u,v$ and $w$ be some PBW monomials in $R$. 
Then
\begin{align}\label{x-bracket}
[ux, wy]&=w[u, y]x+uwx+u[x, w]y,\\
[ux, vxy]&=[uv, x]xy+uvx^2+v[u, y]x^2+v[u, x]x+v[u, y, x]x\nonumber\\
[wy, vxy]&= \alpha vwxy +v[w, x]y^2+v[w, y,x]y\nonumber\\
[v_1xy, v_2xy]&=[v_1v_2, x]zx+[v_1v_2, y]x^2y+(v_1[y, v_2, x]+ v_2[y, v_1, x])xy,\nonumber
\end{align} 
where $z\in \fz$ and $\alpha=1$  if $\deg_H(vw)\in \bZ_0$ and $\alpha=0$ otherwise.
Let $u$ be a PBW monomial in $R$.
Then  $[x,[u,x]]=[R,[u,x]]=0$.
Moreover, if $\deg_H u\in \bZ_0$ then by a simple induction on $\deg_H(u)$ we deduce  that $[u,x]^2=0$. From this it is easy to see that the
associative ideal  $J$ of $u(L)$ generated by all the elements $[u,x]$ 
(with $\deg_H (u)$ even) is nilpotent.
 Consider $u(L)$ as a Lie algebra and set $\fg=u(L)$. 
Since $J$ is a solvable ideal of $\fg$,  it is enough to prove that $\fg/J$ is  solvable.
Now, for all monomials $v,w\in R$, we have
\begin{align}\label{x-bracket-2}
[ux, wy]&=w[u, y]x+uwx+u[x, w]y,\\
[ux, vxy]&=[uv, x]xy \mbox{ modulo } R+J, \nonumber\\
[wy, vxy]&= \alpha vwxy  \mbox{ modulo } R+J, \nonumber\\
[v_1xy, v_2xy]&=[v_1v_2, x]zx+[v_1v_2, y]x^2y+(v_1[y, v_2, x]+ v_2[y, v_1, x])xy,\nonumber
\end{align} 
where $z\in \fz$ and $\alpha=1$  if $\deg_H(vw)\in \bZ_0$ and $\alpha=0$ otherwise.
 It follows from 
Equation \eqref{x-bracket-2} that $[\fg, \fg]\su \mathfrak{m} +J$, where 
$$
\mathfrak{m}= \la ux, wy, vxy, z\mid u, v, w, z\in R,  \deg_H(v)\in \bZ_0 \ra_{\F}.
 $$
Let 
$$
\fk= \la [z_1, v_1xy], v_2xy, z_2\mid v_1, v_2, z_1, z_2\in R,  \deg_H(v_1), \deg_H(v_2) \in \bZ_0 \ra_{\F}.
$$
Now, by using Equation \eqref{x-bracket} one can observe   that $[\fk, \fm]\su \fk+J$ and 
$$\fk^\prime\su \la [z, vxy] \mid  v,  z\in R,   \deg_H(v) \in \bZ_0 \ra_{\F}+J,\quad  \fk^{\prime \prime}\su J.$$
Hence, $\fk+J$ is  a solvable ideal of $\fm+J$  and it is enough to prove $\fm+J/\fk +J$ is  solvable. But $(\fm+J)/(\fk +J) =(\fh+J)/(\fk+J)$, where 
$$
 \fh= \la u x, wy\mid u,w\in R  \ra_{\F}.
 $$
Note that if $\deg_H(u)\in \bZ_0$ and $\deg_H( w)\in \bZ_1$ then $[ux, wy]\in \fk$. 
We deduce that $\fh^\prime\su \la u x, wy\mid u,w\in R, \deg_H(u)\in \bZ_0, \deg_H(w)\in \bZ_1  \ra_{\F}+\fk+J$.
Thus, $\fh^{\prime ^\prime}\su \fk+J$. We deduce that $\fm+J/\fk +J$ is  solvable, as required. 
\qed

It is now a simple matter to prove our main result:

\medspace

\emph{Proof of the Main Theorem.} Note that $u(L)$ is Lie solvable if and only if so is $u(\Lie)$.  The necessary part is then a consequence of Lemma \ref{finite-3} and Lemma \ref{infinite-4}. For the sufficiency observe that if $I$ is a finite-dimensional 2-nilpotent restricted ideal of $\Lie$ then $u(\Lie)$ is Lie solvable if and only if so is 
$u(\Lie)/Iu(\Lie)\cong u(\Lie/I)$. Therefore, if one of the conditions (i)--(v) of the statement holds, then $u(L)$ is Lie solvable by Lemmas \ref{cond-i}, \ref{cond-ii}, \ref{cond-iii}, and \ref{cond-iv}. The proof is complete. \qed

Note that if $C$ is a commutative algebra over a field of characteristic 2 then $R=M_2(C)$ is Lie center-by-metabelian, that
is $[[[R, R], [R, R]], R]=0$. In the  following example we show that  in the statement of the Main Theorem the extension of the ground field is really required:

\begin{example}\label{examplecod} \emph{Let $\F$ be a field of characteristic 2 containing two elements  $\alpha, \beta$ such that the following condition holds: If $\lambda_1,\lambda_2,\lambda_3$ are  in $\F$ and  $\lambda_1^2+\lambda_2^2 \alpha + \lambda_3^2 \beta=0$ then $\lambda_1=\lambda_2=\lambda_3=0$. For instance, one can consider the field $\mathbb{K}(X,Y)$ of rational  functions in two indeterminates over any field $\mathbb{K}$ of characteristic 2, and $\alpha=X$ and $\beta=Y$.  
Let $L$ be the $\F$-vector space  having the elements $x, x_1, x_2, x_3, z_1, z_2, z_3$ as basis . We define a 
 restricted Lie algebra structure on $L$ by setting  
$[x,x_1]=[x, x_3]=z_1$, $[x,x_2]=z_2$,  $ [x_1,x_2]=z_3$, $[x_1, x_3]=\frac{\beta}{\alpha}z_3$, $[x_2,x_3]=0$, $ z_1^{[2]}= z_1, z_2^{[2]}= \alpha  z_1$, $ z_3^{[2]}=\beta z_1$ and $z_i \in Z(L), x_i^{[2]}= 0$, for $i = 1, 2, 3$.\\
\indent Note that $Z(L) = \la  z_1 , z_2, z_3\ra_{\F}$. Furthermore, $Z(L)$ has codimension 4 in $L$ 
and it is free of 2-nilpotent elements. In particular, as $L$ is nilpotent, $Z(L)$ has nontrivial intersection with every nonzero ideal of $L$. As a consequence, $B=0$ is the only 2-nilpotent restricted ideal of $L$, and every 2-abelian restricted ideal of $L$ is indeed abelian. Note also that, as $L$ is nilpotent, none of the condition (iii), (iv), (v) holds for $L$.
We claim that  every abelian restricted ideal of $L$ has codimension at least 2. Indeed, if $L$ has an abelian restricted ideal of codimension 1 then, by  a similar argument as in Lemma \ref{cond-i}, we can embed $u(L)$ 
into $M_2(C)$, for some commutative $\F$-algebra. It follows that $u(L)$
 is Lie center-by-metabelian. However,
one can easily check that $[[x, xx_1], [x_1, x_1x_2x_3], x_2]\neq 0$.\\
\indent Now, let  $\alpha_1$ and $\beta_1$ denote, respectively, the square root of $\alpha$  and $\beta$ in the algebraic closure $\bar{\F}$ of $\F$. Then $v= \alpha_1 z_1+z_2$ and 
$w= \beta_1 z_1+z_3$ are central 2-nilpotent elements of $\Lie=L\otimes_{\F}\bar{\F}$, and 
$J = \bar{\F} v + \bar{\F}w$ is a 2-nilpotent restricted ideal of $\Lie$. Let $I$ be the restricted ideal generated by the images of $x,\alpha x_1+x_2,\beta x_1+x_3$ in $\Lie/J$. Then $I$ is abelian and has codimension 1 in $\Lie/J$, so that condition (i) of the Main  Theorem  assures that $u(L)$ is Lie solvable.}
\end{example}

\section{The ordinary case}
In this section we will use the Main Theorem  in order to characterize Lie solvable universal enveloping algebras of Lie algebras over arbitrary fields.   For this goal the following elementary result is needed. 
Note that, in general, the minimal codimension of abelian ideals of a Lie algebra is not preserved under extensions of the ground field. (See e.g. Example 2.7 of \cite{BC})

\begin{prop}\label{cod1-ordinary} Let $L$ be a Lie algebra over any field $\F$, and let $\bar{\F}$ denote the algebraic closure of $\F$. If $L\otimes_{\F} \bar{\F}$ contains an abelian ideal of codimension 1 then $L$ contains an abelian ideal of codimension 1.  
\end{prop}
\Proof We can obviously assume that $L$ is not abelian. Let ${A}$ be an abelian ideal of $\Lie = L\otimes_{\F} \bar{\F}$ of codimension 1. 
Then $\Lie^\prime \subseteq  A$, in particular $L$ is metabelian. Suppose first that  $C_{\Lie}(\Lie^\prime)= A$ and let $x\in L$ such that $x\otimes 1\notin A$. For every  $y\in L$ there exists $\alpha\in \bar \F$ such that 
$y\otimes 1- x\otimes \alpha\in A$. There exists $z\in L'$ such that $[x\otimes 1, z\otimes 1]\neq 0$.
But $[y\otimes 1- x\otimes \alpha, z\otimes 1]= 0$. We deduce that 
$[y, z]\otimes 1= [x, z]\otimes \alpha$ which implies that $\alpha \in \F$.
Hence, $[y- \alpha x, L^\prime]=0$, which implies that  $C_L(L^\prime)$ is an abelian ideal of $L$ of codimension 1. On the other hand, if $C_{\Lie}(\Lie^\prime)\neq A$ then $C_{\Lie}(\Lie^\prime)=\Lie$, and so $L$ is nilpotent of class 2. If $\dim_{\F}L^\prime=1$ then $\dim_{\bar \F} \Lie^\prime=1$ and it is easy to see that $\dim_{\bar \F} A/Z(\mathfrak{L})=1$. Hence, 
 $\dim_{\F}L/Z(L)=\dim_{\bar \F}\Lie/Z(\Lie)=2$. Now,  the abelian ideal $\la x\ra_{\F} +Z(L)$ has codimension 1 in $L$, for every $x\in L\backslash Z(L)$. Finally, suppose $\dim_{\F}L^\prime>1$. Let $x_1,x_2,y_1,y_2\in L$ such that $z_1=[x_1,y_1]$ and $z_2=[x_2,y_2]$ are linearly independent.  Fix a vector space complement $I_1$  of 
  $\la z_1\ra_{\F}$ in $Z(L)$ such that $z_2\in I_1$. Similarly, we define $I_2$.
  Let $J_i$ be an ideal of $L$ such that $Z(L/I_i)=J_i/I_i$, for $i=1,2$. As $L/I_i$ is not abelian, we have $Z(L/I_i)\otimes_{\F} \bar{\F} \subseteq {A}/(I_i\otimes_{\F }\bar{\F})$. Therefore $J_i$ is abelian. Moreover, since the derived subalgebra of $L/I_i$ is 1-dimensional we have $\dim_{\bar \F}  A / (J_i\otimes_{\F} \bar \F)=1$. Note that, as $z_1$ and $z_2$ are linearly independent modulo $I_1\cap I_2$, one has $J_1\neq J_2$. It follows that $(J_1+J_2)\otimes_{\F} \bar \F= A$, and then $J_1+J_2$ is an abelian ideal of codimension 1 in $L$, yielding the claim. \qed    

\begin{cor} Let $L$ be a Lie algebra over a field ${\mathbb F}$. Then $U(L)$ is Lie solvable if and only if one of the following conditions is satisfied:
\begin{enumerate}
\item [{\normalfont (i)}] $L$ is abelian;
\item [{\normalfont (ii)}] $\char \F=2$ and ${L}$ contains an abelian ideal of codimension 1;
\item [{\normalfont (iii)}] $\char \F=2$, $L$ is nilpotent of class 2 and $\dim_{\F}L/Z(L)=3$;
\item [{\normalfont (iv)}] $\char \F=2$ and $L=\la x_1, x_2, y\ra_{\F}\oplus Z(L)$, with $[x_1,y]=x_1$, $[x_2,y]=x_2$, and $[x_1,x_2]\in Z(L)$.
\end{enumerate}
\end{cor}
\Proof
If $\charac \F \neq 2$ then the assertion is proved in Corollary 6.1 of \cite{RS1}. Then we assume that the ground field has characteristic 2 and $L$ is not abelian. 
Suppose first that $U(L)$ is Lie solvable. Let $\hat{L}$ be the restricted Lie algebra consisting of all primitive elements of the $\F$-Hopf algebra $U(L)$. Then one has  
$$
\hat{L}=\sum_{k \geq 0} {L^{2^k}}\subseteq U(L),
$$
where $L^{2^k}$ is the ${\mathbb F}$-vector space 
spanned by all $x^{2^k}$, where  $x\in L$.
Moreover, $\hat{L}$ is the universal $p$-envelope of $L$ and one has $U(L)=u(\hat{L})$  (see 
 \cite[\S 5.5, Proposition 5.5.3]{M} or  \cite[\S 1.1,  Corollary 1.1.4]{S}), so that the Main Theorem applies. Now, if $\hat{L}\otimes_{\F}\bar{\F}$ contains a 2-abelian ideal of codimension 1, as $u(\hat{L}\otimes_{\F}\bar{\F})=U(L\otimes_{\F}\bar{\F})$ is a domain then it is clear that $L\otimes_{\F}\bar{\F}$ has an abelian ideal of codimension 1. At this stage, Proposition \ref{cod1-ordinary} allows to conclude that $L$ contains an abelian ideal of codimension 1. On the other hand, if condition  (ii) in the statement of the Main Theorem occurs, as $u(\hat{L}\otimes_{\F}\bar{\F})$ is free of nonzero zero divisors we see that $I=0$. Thus $\hat L$ is nilpotent of class 2 and $\dim_{\F}\hat{L}/Z(\hat L)=\dim_{\bar \F}\hat{L}\otimes_{\F} \bar \F/Z(\hat{L}\otimes_{\F} \bar \F)=3$.  It follows at once that $L$ is nilpotent of class 2 and $\dim_{\F}L/Z(L)=3$. Suppose now that condition  (iii) in the statement of the Main Theorem holds. If $L$ has an abelian ideal of codimension 1 then we are done. Otherwise, as $I=0$ and  $\hat{L}\otimes_{\F}\bar \F/Z(\hat L)\otimes_{\F}\bar \F\cong \hat{L}/Z(\hat L)\otimes_{\F}\bar \F$ is 3-dimensional, it is easily seen that $L/Z(L)$  and $\hat{L}/Z(\hat L)$ are isomorphic as ordinary Lie algebras. In particular, as the Lie algebra  $\hat{L}/Z(\hat L)$ is restricted and its derived subalgebra has dimension 2,   by  \cite[Chapter 1, \S 4, (d)]{J} we conclude that $L=\la x_1, x_2, y\ra_{\F}\oplus Z(L)$, with $[x_1,y]=x_1$, $[x_2,y]=x_2$, and $[x_1,x_2]\in Z(L)$. Next note that condition (iv) cannot hold unless $H=0$, which forces  $\dim {\hat L}/Z(\hat L)=2$ and, in turn, $L$ contains an abelian ideal of codimension 1. Finally, suppose that condition (v) in the statement of the Main Theorem  holds.  For every $h\in H$ we have $([x,h]+h)^{[2]}=0$ and so 
$h=[x,h]\in Z(\hat{L}\otimes_{\F} \bar \F)$. Consequently, also in this case we conclude that $L$ contains an abelian ideal of codimension 1. The sufficiency part easily follows from  the Main Theorem and the fact that $u(\hat{L})=U(L)$. \qed

\section*{acknowledgements}
The authors acknowledge the help of Yves Cornulier towards proving Proposition \ref{cod1-ordinary}. They are also grateful to the referee for useful comments.

\end{document}